\documentclass[reqno]{amsart}
\usepackage{amsmath,amsbsy,amsfonts}
\usepackage{color}
\usepackage[utf8]{inputenc}
\usepackage{graphicx}
\usepackage{amsmath,latexsym,amssymb}
\newtheorem{theorem}{Theorem}
\newtheorem{proposition}{Proposition}[section]

\newtheorem{lemma}[proposition]{Lemma}
\newtheorem{remark}[proposition]{Remark}

\newcommand{\diver}{\operatorname{div}}
\newcommand{\dd}{\mathrm{d}}

\numberwithin{equation}{section}

\begin{document}
	
\title[Compact embeddings and fractional spaces]{On compact embeddings into $\mathbf{L^p}$ and fractional spaces}

\author[H.P. Bueno]{H.P. Bueno*}\thanks{* Corresponding author}	
\address[H.P. Bueno]{\newline\indent Departamento de Matemática, Universidade Federal de Minas Gerais, 31270-901 Belo Horizonte, MG, Brazil}
\email{hamilton@mat.ufmg.br}
		
\author[A.H.S. Medeiros]{A.H.S. Medeiros}
\address[A.H.S. Medeiros]{\newline\indent Departamento de Matemática, Universidade Federal de Viçosa, 36570-000 Viçosa, MG, Brazil}
\email{aldo.medeiros@ufv.br}
		
\author[O.H. Miyagaki]{O.H. Miyagaki**}\thanks{** Third author was supported by Grant 2022/16407-1 - São Paulo Research Foundation (FAPESP) and Grant 303256/2022-2 - CNPq/Brazil. }
\address[O.H. Miyagaki]{\newline\indent Departamento de Matemática, Universidade de São Carlos, 13565-905 São Carlos, SP, Brazil}
\email{ohmiyagaki@gmail.com}

\author[G.A. Pereira]{G.A. Pereira}
\address[G.A. Pereira]{\newline\indent Departamento de Matemática, Universidade Federal de Viçosa, 36570-000 Viçosa, MG, Brazil}
\email{gilberto.pereira@ufv.br}


\begin{abstract}
The study of the fractional Laplacian operator $(-\Delta)^s$ in $\mathbb{R}^N$ with Dirichlet boundary conditions gained enormous momentum through its identification with a Neumann operator in $\mathbb{R}^N\times (0, \infty)=\mathbb{R}^{N+1}_+$, a method mainly introduced by Caffarelli and Silvestre. Since then, several other operators have been studied using this method. In general, a crucial question is attached to this method: the  embedding (in the trace sense) on the ground space $L^q(\mathbb{R}^{N})$ is compact? This question is very important  when dealing with problems of existence of solutions. This paper aims to answer this question for some operators.

Passing to an abstract setting, let $X,Y$ be Hilbert spaces and $\mathcal{A}\colon X\to X'$ a continuous and symmetric elliptic operator. We suppose that $X$ is dense in $Y$ and that the embedding $X\subset Y$ is compact. In this paper we show some consequences of this setting for the study of the fractional operator attached to $\mathcal{A}$ in the extension setting $\Omega\times(0,\infty)$ or $\mathbb{R}^{N+1}_+$. Being more specific, we will give some examples where the embedding of the extension domain into  $L^2(\Omega)$ is compact, even in the case $\Omega=\mathbb{R}^N$.	\end{abstract}
	\maketitle
	
\section{Introduction}
The study conducted by Caffarelli and Silvestre \cite{Caffarelli} on fractional powers of the Laplacian operator via a Dirichlet to Neumann operator in the extension setting $\mathbb{R}^N\times (0, \infty)$ was considered by researchers in the field as a breakthrough and led to numerous papers investigating different fractional operators such as $$(-\Delta)^s,\ (-\Delta + m^2)^s,\ \text{and } (-\Delta + |x|^2)^s$$ in the scenario where $0<s<1$. In fact, this study was started by Caffarelli, Salsa and Silvestre in \cite{Salsa}, where it was used to obtain regularity estimates for the obstacle problem in the case of the fractional Laplacian. Stinga and Torrea \cite{StingaTorrea1}, Stinga and Torrea \cite{StingaJFA}, Capella et al. \cite{Capella}, and Arendt et al. \cite{Arendt} have made significant contributions in this area, some of them will be recalled in this paper. Here we address the following question: the embedding of the extension setting $\Omega\times(0,\infty)$ into $L^q(\Omega)$ is compact, even in the case $\Omega=\mathbb{R}^N$?

Let $X$ be a Hilbert space and consider an operator $\mathcal{A}\colon X\to X^{\prime}$, which is supposed to be continuous and symmetric. We admit that $X$ is dense in the Hilbert space $Y$ and that the embedding $X\subset Y$ is compact, see Section \ref{prelim}. Our goal is to show some consequences of this setting on the study of the fractional operator attached to $\mathcal{A}$ in the extension setting $\mathbb{R}^N\times (0, \infty)$.

When dealing with fractional spaces in $\mathbb{R}^N\times (0, \infty)$, the Bessel equation is crucial when defining the space where the fractional operator will be considered:
\begin{equation}\label{Bessel}-u^{\prime\prime}-\frac{1-2s}{t}u^{\prime}+\mathcal{A}u=0,\quad t\in (0,\infty),
\end{equation}

Of course, equation \eqref{Bessel} is related to the extension problem
\begin{equation}\label{Ext}
	\left \{\begin{aligned}
		\displaystyle- \mathcal{A} u +\frac{1-2s}{t}u_{t} + u_{tt}&= 0 \\
		\displaystyle u(x,0) &= f(x).
	\end{aligned}
	\right.
\end{equation}

The well-known case $\mathcal{A}=-\Delta$ is a paradigm. As an example, let us consider it in a bounded Lipschitz domain $\Omega$ in a problem with Dirichlet boundary condition. So, $-\Delta$ is well-defined in $X=H^1_0(\Omega)$, which is compactly embedded into $Y=L^2(\Omega)$. For any $u\in H^1_0(\Omega)$ we have
\[u=\sum_{k=1}^\infty \alpha_k\phi_k\]
where $\phi_k$ are the eigenfunctions of $-\Delta$ in $\Omega$ associated with the eigenvalues $0<\lambda_1<\lambda_2\leq \cdots\leq \lambda_n\to\infty$, the operator $(-\Delta)^s$ can be defined by
\[(-\Delta)^su=\sum_{k=1}^\infty \lambda_k^{s}\alpha_k\phi_k.\]

While studying the problem $(-\Delta)^su=\lambda f(u)$ with Dirichlet boundary conditions, Capella et al. \cite{Capella} considered the space
\[H_{0, L}^1\left(\mathcal{C}_{\Omega}\right)=\left\{v \in L^2\left(\mathcal{C}_{\Omega}\right)\,:\, v=0 \text { on } \partial_L \mathcal{C}_{\Omega}, \iint_{\mathcal{C}_{\Omega}} t^{1-2 \alpha}|\nabla v(x,t)|^2 \dd x \dd t<\infty\right\},\]
(with $\mathcal{C}_{\Omega} = (0,\infty) \times \Omega$ and $\partial_L\mathcal{C}_{\Omega}=\partial\Omega\times [0,\infty)$) and proved that the trace operator $H^1_{0,L}(\mathcal{C}_\Omega)$ to $H^s_0(\Omega)$ is well-defined, see \cite[Section 2.1]{Capella}. 

However, the effect of the compactness of the embedding $H^1_0(\Omega)\subset L^2(\Omega)$ in the embedding $H_{0, L}^1\left(\mathcal{C}_{\Omega}\right)\subset L^q(\Omega)$ for $q\in [2,2^*)$ was not considered and it turns out that the last one is also compact, see Section \ref{Laplace}. As mentioned before, this compactness property might be  fundamental when dealing with existence problems.

In order to obtain this result, we make use of a result obtained by Lions and Magenes \cite[Chap. 1, Théorème 16.2]{LionsMagenes} on the (complex) fractional interpolation space $[X,Y]_{s}$. In fact, they proved that\vspace*{.2cm}

\noindent\textbf{Theorem.} \textit{If the injection $X \rightarrow Y $ is compact, then $$[X,Y]_{\theta_1}\rightarrow [X,Y]_{\theta_2}$$ is also compact, if $0 < \theta_1<\theta_2 < 1 $.}\vspace*{.2cm}

Observe that the compactness of the embedding $H_{0, L}^1\left(\mathcal{C}_{\Omega}\right)\subset L^q(\Omega)$ does not follow immediately from the result obtained by J.-L. Lions and Magenes. This is done in Theorem \ref{Main}, by noting that the space $H_{0, L}^1\left(\mathcal{C}_{\Omega}\right)$ is continuously embedded into the space $\mathcal{H}$ defined there, which, on its turn is compactly embedded into $L^q(\Omega)$, thanks to the result proved by Lions and Magenes.
\vspace*{.3cm}

As a second example, let us consider the fractional powers of the harmonic oscillator, which were considered in Stinga and Torrea \cite[Section 4]{StingaTorrea1} but mainly in the paper by Stinga and Torrea \cite{StingaJFA}. 

The first paper addresses issues related to obtaining fractional operators via the extension setting $\mathbb{R}^{N+1}_+=\mathbb{R}^N\times (0,\infty)$, while \cite{StingaJFA} is devoted to regularity issues concerning the operator $(-\Delta+|x|^2)^s$.

The classical harmonic oscillator $-\Delta+|x|^2$ in $\mathbb{R}^N$ is well-defined in the space $X$ given by
$$X =\left\lbrace u \in H^{1}(\mathbb{R}^{N})\,:\,  \displaystyle\int_{\mathbb{R}^{N}}\left[ |\nabla u|^2 + |x|^2 u^2 \right] dx < \infty \right\rbrace  $$
and we have that the embedding $X\subset Y=L^2(\mathbb{R}^N)$ is compact, a result that we prove here, see Lemma \ref{compact}.

As in \cite{StingaTorrea1} and \cite{StingaJFA}, fractional powers of the harmonic oscillator can be dealt with in the extension setting $\mathbb{R}^{N+1}_+\!=\!\mathbb{R}^N\times (0,\infty)$ and the problem is well-defined in the space
	\[H\!\left(\mathbb{R}^{N+1}_+\right)=\left\{v \in L^2\left(\mathbb{R}^{N+1}_+\right):\iint_{\mathbb{R}^{N+1}_{+}}  \left[\vert \nabla u(x,t) \vert^2 + \vert x \vert^2\vert u(x,t) \vert^2 \right]t^{1-2s} \dd x \dd t<\infty\right\},\]
see Section \ref{harmonic}. As in the previous example, the compactness of the embedding $X\subset Y$ implies the compactness of the embedding $H\left(\mathbb{R}^{N+1}_+\right)\subset L^{q}(\mathbb{R}^{N})$ for every $q \in [2,2^*)$, a result that was not obtained in the aforementioned papers.

Observe that, when dealing with operators in the extension setting $\mathbb{R}^{N+1}_+$, the embedding in the trace sense $H(\mathbb{R}^{N+1}_+)\subset L^2(\mathbb{R}^N)$ (where the space $H(\mathbb{R}^{N+1}_+)$ depends on the operator $\mathcal{A}$) plays an essential role in many problems, since the problem considered has a boundary term in $\mathbb{R}^N$. 

Another examples will be given in Section \ref{consequences}. 

\section{Preliminaries}\label{prelim}

Let  $X,Y$ be Hilbert spaces with $X$ densely embedded into $Y$ and $\mathcal{A}\colon X \to X^{\prime}$ an elliptic operator, which we also suppose to be continuous and symmetric. For any fixed $s \in (0,1)$, let us denote $u(t)=u(t)(x)$ and $A$ the extension of $\mathcal{A}$ to $Y$.

In their influential paper \cite{Arendt},  Arendt et al.  considered problem \eqref{Ext} in the space
\begin{align*}W_s (Y, X):=\Big\{u \in L^{1}_{loc}(X)\, :\, u'\in L^{1}_{loc} (Y),\ t^s u(t) \in L^{*}_{2} (X)\text{ and }t^s u'(t)\in L^*_2 (Y)
\Big\},\end{align*}
where, for $Z=X$ or $Z=Y$
$$L^*_2 (Z):= L_2
\left(Z,\frac{\dd t}{t}\right)=L^2\left((0,\infty);Z,\frac{\dd t}{t}\right).$$
(Observe that our notation differs from \cite{Arendt}, since we have  $X$ densely embedded into $Y$.)

Endowed with the norm
\begin{align*}
	\|u\|_{W^s(Y,X)}&=\left(\|t^su\|^2_{L^*_2(X)}+\|t^su'\|^2_{L^*_2(Y)}\right)^{1/2}\\
	&=\left(\int_0^\infty \Big(\langle u(t),u(t)\rangle_X+\langle u'(t),u'(t)\rangle_Y\Big)t^{2s-1}\dd t\right)^{1/2},
\end{align*}
$W_s (Y, X)$ is in fact a Hilbert space. Clearly $C^\infty_c ([0,\infty); X)$ is a subspace of $W_s (Y, X)$.

A function $u\colon (0,\infty)\to X$ is $s$-\textbf{harmonic} if $u\in W_{1-s} (Y, X)$, $t^{1-2s}u \in W_s (X^{\prime},Y)$ and $u$ is a solution of \eqref{Bessel}, that is
$$-(t^{1-2s} u')' (t) + \frac{1-2s}{t} \mathcal{A}u(t) = 0\quad \text{in } X' \text{ for a.e. }t \in (0, \infty).$$

For each $0<s<1$, $u$ has a normal $s$-derivative given by
$$u^{\prime}(0):=-\lim _{t \downarrow 0} t^{1-2 s} u^{\prime}(t),$$
where the convergence is weak in $X^{\prime}$. It results that $u^{\prime}(0)\in [Y,X^{\prime}]_s$.

In the paper \cite{Arendt}, Arendt et al it is proved that, for any  $z\in [Y, X]_s$, there exists a unique $s$-harmonic function $u$ such that $u(0) = z$.

The authors also considered the Neumann problem, that is, given $y\in [Y,X']_s$, they proved the existence of a unique $s$-harmonic function $u$ such that $$\lim_{t\downarrow 0} -t^{1-2s} u' (t) = y.$$

As a consequence of the well-posedness of these problems, they define a \emph{Dirichlet-to-Neumann} operator $D_s$ in $Y$ satisfying $D_s=c_sA^s$, where $A$ is the extension of $\mathcal{A}$ to $Y$. As it is shown in \cite{Arendt}, the extension obtained by the operator $D_s$ coincides with that of Caffarelli-Silvestre. \vspace*{.2cm}

Let us now suppose that the embedding $X\subset Y$ is also compact, with $X$ dense in $Y$. (The compactness of the embedding $X\subset Y$ was not handled in \cite{Arendt}.) As a consequence, since $\mathcal{A}^{-1}\colon X'\to X$ is a symmetric compact operator, we can assume that the solution $u(x,t)=u(t)(x)$ to \eqref{Ext} has the form
\begin{equation}\label{soma1}
	u(x,t) = \sum_{k=1}^\infty \alpha_k(t) \phi_{k}(x),
\end{equation}
where, for each $k\in \mathbb{N}$, $\alpha_k$ are functions depending only on $t$ and $\{\phi_k\}_{k \in \mathbb{N}}$ is an orthonormal basis of $Y$ formed by autofunctions of $\mathcal{A}$ associated with its eigenvalues
$\lambda_1 < \lambda_2 \leq \cdots\leq \lambda_k \leq \cdots$.

\begin{lemma}\label{l2}
Let us consider the extension problem \eqref{Ext} and suppose that the solution $u$ has the form \eqref{soma1}. Then it holds
\begin{equation*}
\int_{0}^{\infty} \left[ \langle \mathcal{A} u,u \rangle_{Y} + \langle u_t,u_t \rangle_Y \right]t^{1-2s}dt  = K(s)  \sum_{k=1}^\infty \lambda^s_k \langle f, \phi_k\rangle_Y^2,\end{equation*}
where $ K(s) = \frac{-2s\Gamma(-s)}{4^s \Gamma(s)}>0$.
\end{lemma}
\noindent\begin{proof}
In order to simplify our presentation, we return to our previous notation denoting $u(x,t)$ simply by $u(t)$. So,
\[\mathcal{A}u(t) = \sum_{k=1}^\infty \alpha_k(t) \mathcal{A}(\phi_{k})(x)= \sum_{k=1}^\infty \lambda_{k}\alpha_k(t) \phi_{k}(x).
\]
Therefore,
\begin{align*}
	0&= u_{tt} + \frac{1-2s}{t} u_{t} - \mathcal{A}u
	= \sum_{k=1}^\infty \left(\alpha^{\prime\prime}_{k}(t) + \frac{1-2s}{t} \alpha^{\prime}_k(t) + \lambda_{k} \alpha_k(t) \right) \phi_{k}(x)
\end{align*}
and we conclude that for any $t>0$ it holds
\begin{equation}\label{edo}
	-\lambda_{k} \alpha_k+ \frac{1-2s}{t}\alpha^{\prime}_{k} + \alpha^{\prime\prime}_{k}=0.
\end{equation}

We have
\begin{align*}\langle \mathcal{A} u, u \rangle_{Y}
&=\left\langle \sum_{k=1}^\infty \lambda_k \alpha_k(t)  \phi_{k}(x), \sum_{\ell=1}^\infty \alpha_{\ell}(t)  \phi_{\ell}(x)\right\rangle_Y
= \sum_{k=1}^\infty \alpha^2_{k}(t) \lambda_k
\end{align*}
and also
\[\langle u_t, u_t \rangle_{Y} = \sum_{k=1}^\infty \left(\alpha^{\prime}_{k}(t)\right)^2.
\]

Substituting $u(x,0) = f(x)$ into \eqref{soma1}, we obtain
$$\alpha_k(0)= \langle f,\phi_{k} \rangle_{Y}$$
and it follows from \eqref{edo} that
\begin{align*}
	\left(\alpha_k(t)\alpha^{\prime}_{k}(t) t^{1-2s}\right)^{\prime} &= \alpha^{\prime}_{k}(t)  \alpha^{\prime}_{k}(t) t^{1-2s} + \alpha_k(t)  \alpha^{\prime\prime}_{k}(t) t^{1-2s} +(1-2s) \alpha_k(t)  \alpha^{\prime}_{k}(t) t^{-2s}\\
	&=t^{1-2s}\left[ \left(\frac{1-2s}{t} \alpha^{\prime}_{k}(t) + \alpha^{\prime\prime}_{k}(t) \right)\alpha_k(t) + (\alpha^{\prime}_{k}(t))^2\right] \\
	&= t^{1-2s}\left[ \left(\lambda_{k} \alpha_k(t)\right)\alpha_k(t) + (\alpha^{\prime}_{k}(t))^2\right] \\
	&= t^{1-2s}\left[ \lambda_{k} \alpha^2_{k}(t) + (\alpha^{\prime}_{k}(t))^2\right].
\end{align*}

Therefore,
\begin{align*}
	\displaystyle\int_{0}^{\infty} \left[ \langle \mathcal{A} u,u \rangle_{Y} + \langle u_t,u_t \rangle_Y \right]t^{1-2s}dt   &= \displaystyle\int_{0}^{\infty} \left[ \sum_{k=1}^\infty \alpha^2_{k}(t) \lambda_k +  \sum_{k=1}^\infty \left(\alpha^{\prime}_{k}(t)\right)^2 \right]t^{1-2s}dt   \nonumber\\
	&= \sum_{k=1}^\infty \displaystyle\lim_{t \to 0^{+}}\left( -\alpha_k(t)  \alpha^{\prime}_{k}(t) t^{1-2s}\right) \nonumber\\
	&=  \sum_{k=1}^\infty  \frac{-2s\Gamma(-s)}{4^s \Gamma(s)}\lambda^s_k \langle f, \phi_k\rangle_Y^2 \nonumber\\
	&= K(s)  \sum_{k=1}^\infty \lambda^s_k \langle f, \phi_k\rangle_Y^2,
\end{align*}
where $ K(s) = \frac{-2s\Gamma(-s)}{4^s \Gamma(s)}>0$. We observe that the third equality above makes use of $$\lim_{t \to 0^{+}}  \alpha^{\prime}_{k}(t) t^{1-2s} = \frac{-2s\Gamma(-s)}{4^s \Gamma(s)}\lambda^s_k \langle f, \phi_k\rangle,$$
a result which was proved at the end of Section 3.1 in Stinga and Torrea \cite{StingaTorrea1}.
\qed\end{proof}
\section{Main result}
As a consequence of our previous section, we obtain our main result.

\begin{theorem}\label{Main}
Let $X \subset Y$  be Hilbert spaces, the embedding $X\subset Y$ being compact with $X$ dense in $Y$. If $\mathcal{A}:X \to X'$ is an elliptic, continuous and symmetric operator, then the embedding
\[\mathcal{H} \subset Y \]
is compact (in the trace sense), where
$$\mathcal{H}= \left\lbrace u\colon[0,\infty) \to [X,Y]_{s}\,:\, \ \int_{0}^{\infty} \left[\langle u'(t), u'(t) \rangle_Y + \langle\mathcal{A}u(t),u(t) \rangle_Y\right] t^{1-2s}\dd t < \infty \right\rbrace  .$$
\end{theorem}
\noindent\begin{proof} We consider the spaces 
$$X= \left\{ u\, :\, u = \sum_{k=1}^\infty \alpha_{k}\phi_{k}(x), \quad \sum_{k=1}^\infty\lambda^2_k \vert\alpha_{k} \vert^2 < \infty \right\} $$
and
$$[X,Y]_{1- \frac{s}{2}}=  \left\lbrace u \,:\, u = \sum_{k=1}^\infty \alpha_{k} \phi_{k}(x), \quad \sum_{k=1}^\infty\lambda^{s}_k \vert \alpha_{k} \vert^2 < \infty \right\rbrace.$$

It follows from Lemma \ref{l2} that the embedding
$$\mathcal{H} \subset [X,Y]_{1-\frac{s}{2}}$$
is continuous.

Since the embedding $X \subset Y $ is compact, it follows from the Theorem of Lions and Magenes that the embedding  $[X,Y]_{1-\frac{s}{2}}\subset [X,Y]_{1-\frac{s}{4}}$ is also compact.

Now observe that, if $u \in [X,Y]_{1-\frac{s}{4}}$, then
\[\| u \|^2_{[X,Y]_{1-\frac{s}{4}}} = \sum_{k=1}^\infty\lambda^{s}_k |\alpha_k|^2 \geq \lambda^{s}_1 \sum_{k=1}^\infty |\alpha_k|^2 = \lambda_1^s \|u\|^2_{Y},
\]
since 
$\| u \|^2_Y = \sum_{k=1}^\infty |\alpha_k|^2$, proving that $ [X,Y]_{1-\frac{s}{4}} \subset Y$ is continuous.

Therefore, we have the embeddings
$$\mathcal{H} \stackrel{continuous}{\subset} [X,Y]_{1-\frac{s}{2}} \stackrel{compact}{\subset}[X,Y]_{1-\frac{s}{4}} \stackrel{continuous}{\subset}Y,$$
proving that the embedding
$$\mathcal{H} \subset Y$$
is compact in the trace sense. We are done.
\qed
\end{proof}
\section{Consequences}\label{consequences}
\subsection{The Laplacian}\label{Laplace}
Let $\Omega \subset \mathbb{R}^{N}$ be a bounded domain with Lipschitz boundary and $\mathcal{A} = -\Delta$ the Laplacian operator in a problem with Dirichlet boundary condition.  The operator $\mathcal{A}=-\Delta$ is elliptic, symmetric and continuous and it is well-defined in $X= H_0^{1}(\Omega) \subset L^2(\Omega) = Y$, the embedding $X\subset Y$ being compact. This section demonstrates the practical application of our results in this specific situation.

Of course, the literature dedicated to problems involving the fractional Laplacian operator $(-\Delta)^s$ is so broad that we will not even try to quote some of the main contributions to the area. We will restrict ourselves to the study made by Capella et al. \cite{Capella}, which considers the case $\Omega=B_1=B_1(0)\subset\mathbb{R}^N$ with Dirichlet boundary conditions and the problem $(-\Delta)^s=\lambda f(u)$. Adapting their notation to ours, their analysis is carried out in the space (see \cite[Eq. (1.3)]{Capella})
\[H=\left\{u\in L^2(B_1)\,:\, \|u\|_{\mathcal H}=\sum_{k=1}^{\infty}\lambda_k^s|\alpha_k|^2<\infty\right\},\]
which coincides with the space $[X,Y]_{1-s/2}$ defined in the previous section in the case $\Omega=B_1$. 

Observe that, maintaining the notation of our theorem, if $u\in Y$, then
\[
\langle \mathcal{A}u,u \rangle_Y = \int_{\Omega} \left(-\Delta_x u\right)u \dd x = \int_{\Omega} \vert \nabla_xu(x,t) \vert^2 \dd x
\]
and
\begin{multline*}
\int_{0}^{\infty} \left[\langle u'(t), u'(t) \rangle_Y + \langle\mathcal{A}u(t),u(t) \rangle_Y\right] t^{1-2s}\dd t
\end{multline*}
\begin{align*}
&= \int_{0}^{\infty}  \int_{\Omega} \left[\vert u_t(x,t) \vert^2 \dd x + \vert \nabla_x u(x,t) \vert^2 \right]t^{1-2s} \dd x \dd t\nonumber\\
&= \iint_{(0,\infty) \times \Omega}  \vert \nabla u(x,t) \vert^2 t^{1-2s} \dd x \dd t.
 \end{align*}

Thus, the space $\mathcal{H}$ of our theorem is given by
$$\mathcal{H} = \left\{u:(0,\infty) \to [X,Y]_{1-\frac{s}{2}}\,:\, \ \iint_{\mathcal{C}_{\Omega}}  \vert \nabla u(x,t) \vert^2 t^{1-2s} \dd x \dd t < \infty\right\},$$
where $\mathcal{C}_{\Omega} = (0,\infty) \times \Omega$. We denote by $\partial_L\mathcal{C}_{\Omega}$ the boundary of the cylinder $\mathcal{C}_{\Omega}$.

As in Cappela et al. (see p. 1357), if we denote
$$H_{0, L}^1\left(\mathcal{C}_{\Omega}\right)=\left\{v \in L^2\left(\mathcal{C}_{\Omega}\right)\,:\, v=0 \text { on } \partial_L \mathcal{C}_{\Omega}, \iint_{\mathcal{C}_{\Omega}} t^{1-2 \alpha}|\nabla v(x,t)|^2 \dd x \dd t<\infty\right\},
$$
then $H_{0, L}^1\left(\mathcal{C}_{\Omega}\right) \subset \mathcal{H}$ and consequently our result allows us to conclude that the embedding (in the trace sense)
\[H_{0, L}^1\left(\mathcal{C}_{\Omega}\right) \subset L^{2}(\Omega)\]
is compact.

We also know that $H_{0, L}^1\left(\mathcal{C}_{\Omega}\right) \subset L^{2^*}(\Omega)$ is continuous and we conclude that $$H_{0, L}^1\left(\mathcal{C}_{\Omega}\right) \subset L^{q}(\Omega)$$
is compact for each $q \in [2,2^*)$, a result that complements Proposition 2.1 in Cappela et al. \cite{Capella}.

\subsection{Magnetic Operator}
Let $A:\overline{\Omega} \rightarrow \mathbb{R}^N$ be a continuous magnetic potential on a bounded smooth domain $\Omega \subset \mathbb{R}^N$, $N\geq 3$. We consider the Schrödinger operator $L_A=(-i \nabla -A(x))^2$, which is defined in $\Omega$ by
$$L_A u=(-i \nabla -A(x))^2 u = -\Delta u +|A|^2 u +2i A\cdot  \nabla u +iu \diver  A.$$

Formally, the quadratic form associated with the self-adjoint operator $L_A$ is given by
\begin{equation*}
		\langle L_A u, u\rangle = \int_{\Omega}|-i\nabla u -A u|^2 dx=  \int_{\Omega}|\nabla_A u|^2 \dd x.
\end{equation*}
We shall denote by  $H_A (\Omega)$   the Hilbert space obtained as the completion of $C_{0}^{\infty}(\Omega,\mathbb{C})$ under the norm $$\|u\|_A=\left(\int_{\Omega}|\nabla_A u|^2 \dd x\right)^{1/2}.$$

The diamagnetic inequality, proved by Esteban and Lions 	in \cite[Section II]{Esteban}, guarantees that for any $ u \in H_A (\Omega)$ holds
\begin{align}\label{diamagnetic}
|\nabla|u|(x)|&=\left|\mathfrak{Re}\left(\nabla u \frac{\overline{u}}{u}\right)\right|=\left|\mathfrak{Re}\left((\nabla u-i Au)\frac{\overline{u}}{u}\right)\right| \nonumber\\  &\leq  |\nabla_A u(x)|,
\end{align}
where  $\mathfrak{Re} (z)$ and $\overline{z}$ denote the real part and the complex conjugated of $z\in \mathbb{C},$ respectively.
	
If $ u \in H_A (\Omega)$, it follows from \eqref{diamagnetic} that $|u|\in H^{1}_{0}(\Omega)$. Consequently, if $N\geq 3$, the embedding
$$ H_A(\Omega) \subset L^q(\Omega,\mathbb{C}),$$
is continuous for $1\leq q\leq 2^*=\frac{2N}{N-2}$ and compact for $1\leq q< 2^*$. Here $$L^q(\Omega,\mathbb{C})=\left\{u\colon \Omega \rightarrow \mathbb{C}\,:\, \int_{\Omega}|u|^q \dd x< \infty\right\}.$$

With respect to the literature considering a magnetic operator, we would like to cite  \cite{AlvesGio,Arioli,Barile,Chacha,Ding} for problems involving the existence and other properties of solutions with the nonlinearity having polynomial growth, and mention \cite{BarileF} for the case of exponential growth of the nonlinearity. In the case of the fractional magnetic operator, we cite for instance  \cite{Davenia,Manasses} when the nonlinearity has critical polynomial and exponential growth, respectively. In \cite{Wang} is treated a fractional magnetic problem involving Choquard and polynomial perturbation in the Hardy-Littlewood-Sobolev sense.

For all $u \in Y=L^2(\Omega,\mathbb{C})$ we have
$$\langle \mathcal{A}u,u \rangle_{Y}= \langle L_Au, u \rangle = \int_{\Omega} \vert \nabla_A u(x,t) \vert^2  \dd x$$
and so
\begin{align*}
		\lefteqn{\int_{0}^{\infty} \left[\langle u'(t), u'(t) \rangle_Y + \ \langle\mathcal{A}u(t),u(t) \rangle_Y\right] t^{1-2s}\dd t} \\ &=\int_{0}^{\infty}  \int_{\Omega} \left[\vert u_t(x,t) \vert^2 \dd x + \vert \nabla_A u(x,t) \vert^2 +  \vert u(x,t) \vert^2 \right]t^{1-2s} \dd z \dd t \equiv I_{ut}.
\end{align*}

Therefore, the space $\mathcal{H}$ of our theorem is given by
$$\mathcal{H}= \left\{ u\colon[0,\infty) \to [X,Y]_{s}\,:\, I_{ut}< \infty \right\}.$$

Maintaining the notation introduced in the previous example and denoting
$$H\left(\mathcal{C}_{\Omega}\right)=\left\{u \in L^2\left(\mathcal{C}_{\Omega}\right)\,:\,I_{ut}<\infty\right\},
$$
we have the embedding $H\left(\mathcal{C}_{\Omega}\right) \subset \mathcal{H}$ and we conclude that, in the trace sense, the embedding
$$H\left(\mathcal{C}_{\Omega}\right) \subset L^{q}(\Omega)$$
is compact. Since the embedding $H\left(\mathcal{C}_{\Omega}\right) \subset L^{2^{*}_\gamma}(\Omega)$ is continuous, we conclude that the embedding $H\left(\mathcal{C}_{\Omega}\right) \subset L^{q}(\Omega)$ is compact for every $q \in [2,2^{*})$.

\subsection{Grushin Operator}

Let $\Omega\subset \mathbb{R}^{m+n}=\mathbb{R}^m\times\mathbb{R}^n$ be a smooth bounded domain and let us  denote $z=(x,y)\in \mathbb{R}^m\times\mathbb{R}^n$.

The Grushin operator $\mathcal{A}=-\Delta_\gamma $  is given by $\Delta_\gamma u:= \Delta_x u +(1+\gamma)^2 |x|^{2 \gamma}\Delta_y u$, where $\Delta_x$ and $\Delta_y$  are the Laplacian operators in the variables $x$ and $y$, respectively, and  $\gamma >0$. The Grushin operator is well-defined in the space $X:=\dot{H}^{1,2}_{\gamma}(\Omega)$,  which is the completion of $C_{0}^{\infty}(\Omega)$ with respect to the norm
$$\|u\|:=\left(\int_\Omega |\nabla_\gamma u|^2 \dd z\right)^{1/2},$$
where $\nabla_\gamma u=(\nabla_x u, (1+\gamma)|x|^{\gamma} \nabla_y u)$.

We recall that this operator is elliptic for $ x \neq  0$  and degenerates on the manifold $\{0\}\times \mathbb{R}^n$. When $\gamma$ is a nonnegative integer, this operator belongs to the class of Hörmander-type operators. In the general case, the operator belongs to a class of sub-elliptic operators studied by Franchi and Lanconelli in \cite{[11], [12],[13]}. Alternatively, it can be seen as an $X$-elliptic  operator, which was treated  by Kogoj and Lanconelli in \cite{Kogoj} and \cite{[19]}, respectively. Louidice  \cite{Louidice} studied  semilinear sub-elliptic equations with Hardy term and critical nonlinearity, and Alves and Holanda \cite{AH} considered a zero mass problem.

Keeping the notation introduced before, for $N=m+n$, the embedding $X=\dot{H}^{1,2}_{\gamma}(\Omega)\subset Y=L^2\left(\mathcal{C}_{\Omega}\right)$ is compact, see \cite[Prop 3.2]{Kogoj} and also \cite[Prop 2.1]{AlvesKogoj}.

For all $u \in Y$ we have
\begin{align*}
\langle \mathcal{A}u,u \rangle_{Y} &= \int_{\Omega} ( \Delta_x  +(1+\gamma)^2 |x|^{2 \gamma}\Delta_y )u^2 \dd z\\ & = \int_{\Omega} \vert \nabla_x u(x,t) \vert^2 + (1+\gamma)^2 |x|^{2 \gamma} \vert \nabla_y  u(x,t) \vert^2 \dd z
\end{align*}
what yields
\begin{align*}
		\lefteqn{\int_{0}^{\infty} \left[\langle u'(t), u'(t) \rangle_Y + \ \langle\mathcal{A}u(t),u(t) \rangle_Y\right] t^{1-2s}\dd t} \\ &=\int_{0}^{\infty}  \int_{\Omega} \left[\vert u_t(x,t) \vert^2 \dd x + \vert \nabla_x u(x,t) \vert^2 +  (1+\gamma)^2 |x|^{2 \gamma} \vert u(x,t) \vert^2 \right]t^{1-2s} \dd z \dd t \\
		&= \iint_{\mathcal{C}_{\Omega}}  \left[\vert \nabla u(x,t) \vert^2 +  (1+\gamma)^2 |x|^{2 \gamma}\vert u(x,t) \vert^2 \right]t^{1-2s} \dd z \dd t\equiv I_{ut}
\end{align*} 

So, the space $\mathcal{H}$ of our theorem is given by $$\mathcal{H}= \left\{ u\colon[0,\infty) \to [X,Y]_{s}\,:\, I_{ut}< \infty \right\}.$$
	
Denoting
$$H\left(\mathcal{C}_{\Omega}\right)=\left\{u \in L^2\left(\mathcal{C}_{\Omega}\right):I_{ut}<\infty\right\},$$
we have the embedding $H\left(\mathcal{C}_{\Omega}\right) \subset \mathcal{H}$ and we conclude that, in the trace sense, the embedding
$$H\left(\mathcal{C}_{\Omega}\right) \subset L^{q}(\Omega)$$
is compact.

Since the embedding $H\left(\mathcal{C}_{\Omega}\right) \subset L^{2^{*}_\gamma}(\Omega)$ is continuous, we conclude that the embedding $H\left(\mathcal{C}_{\Omega}\right) \subset L^{q}(\Omega)$ is compact for every $q \in [2,2^{*}_\gamma)$, where $2^{*}_\gamma=\frac{2N_\gamma}{N_\gamma-2}$ and $N_\gamma=m+(1+\gamma)n.$

\subsection{The harmonic oscillator}\label{harmonic}

We now consider the operator $\mathcal{A} = -\Delta + \vert x \vert^2$ in $\Omega=\mathbb{R}^N$. Fractional powers of the harmonic oscillator, namely $(-\Delta+|x|^2)^s$, $0<s<1$, have also been considered. Applying mathematical physics methods, the special case $s=1/2$ was handled by Delbourgo \cite{Delbourgo}, but his analysis is quite different from our approach.

The general case $0<s<1$ was considered by Stinga and Torrea in the influential paper \cite{StingaTorrea1}, where among other notable results, a Poisson formula for general self-adjoint lower bounded operators and a Harnack inequality for the fractional harmonic oscillator were proved.

The study of the operator $(-\Delta+|x|^2)^s$ was further developed in another important paper by Stinga and Torrea \cite{StingaJFA} by defining a different class of Hölder spaces $C_H^{k,\alpha}$, which preserves the Hermite–Riesz transforms.  Results concerning both Schauder estimates and the regularity of $$(-\Delta+|x|^2)^s u$$  were obtained, and also a pointwise definition of the operator for $u\in C_H^{k,\alpha}$. These papers are sequels to a previous paper by Bongioanni and Torrea \cite{BongioanniTorrea}, where fractional inverses $H^{-a}$ of the operator $H=-\Delta+|x|^2$ where studied for $a>0$.

Despite having different characteristics, the operator $\left(-\Delta +|x|^2\right)^s$ clearly recalls the operator $\left(-\Delta+m^2\right)^s$, which has been the subject of a number of papers in recent years, starting with the case $s=1/2$ in the paper by Coti Zelati and Nolasco \cite{ZelatiNolasco}, which was generalized in \cite{Cingolani}, see also \cite{BBMP}. The operator $\sqrt{-\Delta+m^2}$ is known as the \emph{pseudo-relativistic Hamiltonian operator} and describes, in the mathematical context, the Hamiltonian for the motion of a free relativistic particle and therefore has broad applications in Physics. A generalization for any $0<s<1$, that is  $(-\Delta+m^2)^s$, is also the subject of various papers: see, e.g. \cite{Ambrosio,ABMP, BOMP,BMP,Fall}. 

\begin{remark}
Observe that, considering the operator $(-\Delta+m^2)^s$ $(0<s<1)$ in a smooth, bounded domain $\Omega\subset\mathbb{R}^N$, the same arguments of Section \ref{Laplace} can be carried out.
\end{remark}

The operator $\mathcal{A}=-\Delta +|x|^2$ is well-defined in the space $X$ given by
$$X =\left\lbrace u \in H^{1}(\mathbb{R}^{N})\,:\,  \displaystyle\int_{\mathbb{R}^{N}}\left[ |\nabla u|^2 + |x|^2 u^2 \right] \dd x < \infty \right\rbrace .$$

\begin{lemma}\label{compact} The embedding $X \hookrightarrow L^{2}(\mathbb{R}^N)$ is continuous for all $q\in [2, 2^*]$ and compact for any $q\in [2, 2^*)$.
\end{lemma}
\noindent\begin{proof} 
For any $u \in X$ we have
\begin{align*}
\int_{\mathbb{R}^N}|u|^2 \dd x &= \int_{|x| \leq 1}|u|^2 \dd x \ + \ \int_{|x|>1}|u|^2 \dd x \leq \int_{|x| \leq 1}|u|^2 \dd x + \int_{|x|>1} |x|^2|u|^2 \dd x \\
&\leq \int_{|x| \leq 1}|\nabla u|^2 \dd x + \int_{|x|>1} |x|^2|u|^2 \dd x \\
&\leq K \int_{ \mathbb{R}^{N}} \left[ |\nabla u|^2 + |x|^2 u^2 \right] \dd x,
\end{align*}
as a consequence of the continuity of the Sobolev embedding for bounded domains. 

On the other hand, the Sobolev-Gagliardo-Nirenberg inequality asserts that there exists a positive constant $S$ such that
$$\int_{\mathbb{R}^N}|u|^{2^*} \dd x \leq S \int_{\mathbb{R}^N}|\nabla u|^2 \dd x \leq S \int_{ \mathbb{R}^{N}} \left[ |\nabla u|^2 + |x|^2 u^2 \right] \dd x.$$
	
Therefore we have proved the continuity of the embedding for $q=2$ and $q=2^*$. The continuity of the embedding for a fixed $q \in\left(2,2^*\right)$ follows from the interpolation inequality, see below.

Let us now assume that $(u_n) \subset X$ satisfies $u_n \rightharpoonup 0$. We will show that  \begin{equation}\label{L2}
	u_n\rightarrow 0\qquad\textrm{in }\ L^{2}(\mathbb{R}^N).
\end{equation}

Since $u_n \rightharpoonup 0$ in $X$, there exists $C$ such that, for all $n\in\mathbb{N}$,
\begin{equation}\label{oli2}
	\int_{ \mathbb{R}^{N}} \left[ |\nabla u_n|^2 + |x|^2 u^2_n \right] \dd x\leq C.
\end{equation}

For any $\varepsilon >0$ fixed, take $R>0$ such that $ \frac{C}{R^2} < \varepsilon$. It follows from \eqref{oli2} that
\begin{align*}
\int_{ \mathbb{R}^N } |u_n(x)|^2 \dd x & =\int_{ B_R (0) } |u_n(x)|^2 \dd x + \int_{ B^c_R (0)} |u_n(x)|^2 \dd x \\
&\leq \int_{ B_R (0) } |u_n(x)|^2 \dd x + \int_{ B^c_R (0)} \frac{|x|^2}{R^2}|u_n(x)|^2 \dd x \\
&\leq \int_{ B_R (0) } |u_n(x)|^2 \dd x + \frac{1}{R^2} \int_{ B^c_R (0)}  \left[|\nabla u_n |^2 + |x|^2|u_n(x)|^2\right] \dd x \\
&\leq \int_{ B_R (0) } |u_n(x)|^2 \dd x +\frac{C}{R^2}\\
&\leq \int_{ B_R (0) } |u_n(x)|^2 \dd x +\varepsilon.
\end{align*}

Since $u_n \rightarrow 0$ in $L^2_{loc}(\mathbb{R}^N)$, passing to the limit as $n\to\infty$ yields 
$$\limsup_{n\in\mathbb{N}}\int_{ \mathbb{R}^N } |u_n(x)|^2 \dd x \leq  \varepsilon.$$

We conclude that $\|u_n\|_{2} \rightarrow 0$, completing the proof of \eqref{L2}.

If $ q \in (2, 2^{*}) $, then 
$$\|u_n\|_{q} \leq \|u_n\|^{1-t}_{2} \|u_n\|^{t}_{2^{*}},$$
where $t \in (0,1)$ and $\frac{1}{q} = \frac{1-t}{2} + \frac{t}{2^{*}}$.

Thus, $\|u_n\|_{q} \rightarrow 0$ and we are done.
$\hfill\Box$\end{proof}\vspace*{.3cm}

Once Lemma \ref{compact} is proved, we return to the study of the operator $(-\Delta+|x|^2)^s$ in $\mathbb{R}^N$ in order to prove the compactness of the embedding into $Y=L^2(\mathbb{R}^N)$. In fact, for all $u \in Y$ we have
\[
\langle \mathcal{A}u,u \rangle_{Y} = \int_{\mathbb{R}^{N}} (-\Delta_x + \vert x \vert^2)u^2 \dd x = \int_{\mathbb{R}^{N}} \vert \nabla_x u(x,t) \vert^2 + \vert x \vert^2 \vert u(x,t) \vert^2 \dd x
\]
and so
\begin{multline*}
	\int_{0}^{\infty} \left[\langle u'(t), u'(t) \rangle_Y+ \ \langle\mathcal{A}u(t),u(t) \rangle_Y\right] t^{1-2s}\dd t 
\end{multline*}
\begin{align*}
&=\int_{0}^{\infty}  \int_{\mathbb{R}^{N}} \left[\vert u_t(x,t) \vert^2 \dd x + \vert \nabla_x u(x,t) \vert^2 + \vert x \vert^2 \vert u(x,t) \vert^2 \right]t^{1-2s} \dd x \dd t\\
&= \iint_{\mathbb{R}^{N+1}_{+}}  \left[\vert \nabla u(x,t) \vert^2 + \vert x \vert^2\vert u(x,t) \vert^2 \right]t^{1-2s} \dd x \dd t.
\end{align*}

So, the space $\mathcal{H}$ of our theorem is given by
\[\mathcal{H}= \left\lbrace u\colon[0,\infty) \to [X,Y]_{s}\,:\, \iint_{\mathbb{R}^{N+1}_{+}}  \left[\vert \nabla u(x,t) \vert^2 + \vert x \vert^2\vert u(x,t) \vert^2 \right]t^{1-2s} \dd x \dd t < \infty \right\rbrace .
\]
Denoting
$$H\!\left(\mathbb{R}^{N+1}_+\right)\!=\!\left\{u \in L^2\left(\mathbb{R}^{N+1}_+\right)\!:\!\iint_{\mathbb{R}^{N+1}_{+}}\!  \left[\vert \nabla u(x,t) \vert^2 + \vert x \vert^2\vert u(x,t) \vert^2 \right]\!t^{1-2s} \dd x \dd t<\infty\right\},
$$
we have the embedding $H\left(\mathbb{R}^{N+1}_+\right) \subset \mathcal{H}$ and we conclude that, in the trace sense, the embedding
\[H\left(\mathbb{R}^{N+1}_+\right) \subset L^{2}(\mathbb{R}^{N})\]
is compact.

Since the embedding $H\left(\mathbb{R}^{N+1}_+\right) \subset L^{2^*}(\mathbb{R}^{N})$ is continuous, we conclude that the embedding $H\left(\mathbb{R}^{N+1}_+\right) \subset L^{q}(\mathbb{R}^{N})$ is compact for every $q \in [2,2^*)$.

\subsection{On a traveling wave equation of the Klein-Gordon type}
Let us consider the operator
\[-\Delta u+b(x)u \quad\text{in }\ \mathbb{R}^N,\]
where we assume that $b(x)$ satisfies
\begin{enumerate}
	\item [$(b_1)$] $b\in C(\mathbb{R}^N,\mathbb{R})$ and $\inf_{x\in\mathbb{R}^N}b(x)=b_0>0$;
	\item [$(b_2)$] for every $M>0$, if $\mu$ denotes de Lebesgue measure in $\mathbb{R}^N$, then
	\[\mu\left(\left\{x\in\mathbb{R}^N\,:\, b(x)\leq M\right\}\right)<\infty.\]
\end{enumerate}
This was one of the operators considered by Bartsch and Wang in \cite{Bartsch}, where the problem
\begin{equation}\label{Bartsch}-\Delta u+b(x)u=f(x,u)\end{equation}
was studied for a continuous, superlinear and subcritical nonlinearity. In that paper, considering the space
\[X=\left\{u\in H^1(\mathbb{R}^N)\,:\, \int_{\mathbb{R}^N}[|\nabla u|^2+b(x)u^2]\dd x<\infty\right\},\]
the authors prove that problem \eqref{Bartsch} has both a positive and a negative weak solution and, if $f(x,-u)=-f(x,u)$ that \eqref{Bartsch} has in fact infinitely many solutions.

In their proof of the above mentioned result, Bartsch and Wang conclude that the embedding $X\subset L^2(\mathbb{R}^N)$ is compact (see \cite[Remark 3.5]{Bartsch}) and, by applying the Sobolev-Gagliardo-Nirenberg inequality, that $X\subset L^p(\mathbb{R}^N)$ is also compact for any $2\leq p<2N/(N-2)$.

Applying the same reasoning of the previous subsection, we examine the operator $(-\Delta+b(x))^s$ in the extension space $\mathbb{R}^{N+1}_+$. Denoting
$$H\!\left(\mathbb{R}^{N+1}_+\right)\!=\left\{u \in L^2\left(\mathbb{R}^{N+1}_+\right):\iint_{\mathbb{R}^{N+1}_{+}}  \left[\vert \nabla u(x,t) \vert^2 + b(x)\vert u(x,t) \vert^2 \right]t^{1-2s} \dd x \dd t<\infty\right\},$$
it turns out, as before, that the embedding  $H\left(\mathbb{R}^{N+1}_+\right)\subset L^2(\mathbb{R}^N)$ is compact.

\subsection{A semilinear Schrödinger equation in $\mathbb{R}^2$}
We recall the following result, which is an adaptation to our setting of a more general outcome obtained by Yunyan Yang in \cite{Yang}.

\noindent\textbf{Theorem.} \textit{Assume that $V\colon \mathbb{R}^2\to \mathbb{R}$ is continuous and satisfies
\begin{enumerate}
	\item [$(V_1)$] $V(x)\geq V_0>0$ for some constant $V_0$;
	\item [$(V_2)$] The function $\displaystyle\frac{1}{V(x)}\in L^1(\mathbb{R}^2)$.
\end{enumerate}
	Then the space
	\[X=\left\{u\in H^1(\mathbb{R}^2)\,:\,\int_{\mathbb{R}^2}V(x)u^2\dd x<\infty\right\}\]
	is compactly embedded into $L^q(\mathbb{R}^2)$ for all $q\geq 1$.}\vspace*{.2cm}

As a consequence of Theorem \ref{Main}, we consider the operator $(-\Delta +V(x))^s$ in the extension space $\mathbb{R}^3_+$. Proceeding as in Subsection \ref{harmonic}, the space $\mathcal{H}$ of our theorem is given by
\[\mathcal{H}=\left\lbrace u\colon[0,\infty) \to [X,Y]_{s}\,:\, \iint_{\mathbb{R}^{3}_{+}}  \left[\vert \nabla u(x,t) \vert^2 + V(x)\vert u(x,t) \vert^2 \right]t^{1-2s} \dd x \dd t < \infty \right\rbrace .\]

Denoting
$$H\left(\mathbb{R}^{3}_+\right)=\left\{u \in L^2\left(\mathbb{R}^{3}_+\right):\iint_{\mathbb{R}^{3}_{+}}  \left[\vert \nabla u(x,t) \vert^2 + V(x)\vert u(x,t) \vert^2 \right]t^{1-2s} \dd x \dd t<\infty\right\}
$$
we have the embedding $H\left(\mathbb{R}^{3}_+\right) \subset \mathcal{H}$ and we conclude that, in the trace sense, the embedding
\[H\left(\mathbb{R}^{3}_+\right) \subset L^{2}(\mathbb{R}^{2})\]
is compact.

In the aforementioned paper, Yang considers the problem
\begin{equation}\label{Yangeq}-\Delta_Nu+V(x)|u|^{N-2}u=\frac{f(x,u)}{|x|^\beta} \text{in }\ \mathbb{R}^N,\end{equation}
where $0<\beta<N$, $V$ satisfies $(V_1)$ and $(V(x))^{-1}\in L^{1/(N-1)}(\mathbb{R}^N)$ and $f$ is continuous and behaves like $e^{\alpha s^{N/(N-1)}}$ as $s\to \infty$. It is proved that \eqref{Yangeq} has a nontrivial positive mountain-pass type weak solution. Observe that the space considered in the paper is
\[E=\left\{u\in W^{1,N}(\mathbb{R}^N)\,:\, \int_{\mathbb{R}^N}V(x)|u|^N\dd x\leq \infty\right\}.\]

\subsection*{Competing interests}
The authors declare that they have no conflict of interests.

\subsection*{Authors' contributions}
All authors contributed equally to the article.
\subsection*{Availability of data and materials}
This declaration is not applicable.

\end{document}